\documentclass{amsart}

\usepackage{graphicx}
\usepackage{textcomp}

\newtheorem{thm}{Theorem}[section]

\newtheorem{prp}{Proposition}[section]

\newtheorem{rem}{Remark}[section]
\newtheorem{dfn}{Definition}[section]

\begin{document}
\title
[Sercan Turhan\ {*}, Nurg\"{u}l OKUR BEKAR, Hande G\"{U}NAY AKDEM\.{I}R]
{Hermite-Hadamard Type Inequality for Log-convex Functions via Sugeno Integrals
}
\author{Sercan Turhan\ {*}}
\address{Dereli Vocational School, Restaurant and Catering Cookery Program,
Giresun University, 28100, Giresun, Turkey.}
\email{sercanturhan28@gmail.com}

\author{Nurg\"{u}l OKUR BEKAR}
\address{Department of Statistics, Faculty of Arts and Sciences,
Giresun University, 28100, Giresun, Turkey.}
\email{nrgokur@gmail.com}

\author{Hande G\"{U}NAY AKDEM\.{I}R}
\address{Department of Mathematics, Faculty of Arts and Sciences,
Giresun University, 28100, Giresun, Turkey.}
\email{hande.akdemir@giresun.edu.tr}

\keywords{Hermite Hadamard type inequality; Sugeno integrals;
log-convex functions.
*Corresponding author}

\begin{abstract}
 In this paper, Hermite-Hadamard type inequality for Sugeno integrals based on log-convex functions
 is studied. Some examples are given to illustrate the results.
\end{abstract}
\maketitle

\section{Introduction}
The theory of fuzzy measures and fuzzy integrals was introduced by
Sugeno \cite{Sugeno}. The properties and applications of
Sugeno-integral have been studied by lots of authors. Between these
others, Ralescu and Adams \cite{Ralescu} proposed several equivalent
definitions of fuzzy integrals; Rom\'{a}n-Flores et al. \cite{Dubois,Flores}
defined the level-continuity of fuzzy integrals and the H-continuity
of fuzzy measures; the book by Wang and Klir \cite{Wang} contains a general
overview on fuzzy measurement and fuzzy integration theory.

Many authors generalized the Sugeno integral by using some other
operators to replace the special operators $\vee$ and/or $\wedge.$
Su\'{a}rez Garc\'{\i}a and Gil \'{A}lvarez \cite{Surez} presented two
families of fuzzy integrals, the so-called seminormed fuzzy
integrals and semiconormed fuzzy integrals. 

In recent years, some authors \cite{Caballero_1}-\cite{Caballero_3} generalized several classical
integral inequalities for fuzzy integral.Caballero and Sadarangani
\cite{Caballero_3} showed off a Hermite-Hadamard type inequality of fuzzy
integrals for convex function. Li, Song and Yue \cite{Li} served
Hermite-Hadamard type inequality for Sugeno integrals. In \cite{Dragomir_1},
Dragomir and Mond introduced to Hermite-Hadamard type inequality for
log-convex functions.

The aim of this paper is to prove a Hermite-Hadamard type inequality
for Sugeno integrals related to log-convex functions. Some example
are given to illustrate the results.

Let's see some proporties of fuzzy integral.
\vskip 0.75cm
\section{Preliminary Discussions}

In this section, we remember some basic definition and properties of
fuzzy integral and log-convex function. For details we refer the
readers to Refs \cite{Sugeno, Wang, Li}.

Suppose that $\Sigma$ is a $\sigma$-algebra of subsets of $X$ and that $\mu:\Sigma\rightarrow\lbrack0,\infty)$ is a
non-negative, extended real-valued set function. We say that $\mu$
is a fuzzy measure if and only if:

\begin{enumerate}
\item $\mu\left(  \emptyset\right)  =0$;

\item $E,F\in\Sigma$ and $E\subset F$ imply $\mu\left(  E\right)  \leq
\mu\left(  F\right)  $ (monotonicity);

\item $\left\{  E_{n}\right\}  \subset\Sigma,E_{1}\subset E_{2}\subset...,$
imply $\underset{n\rightarrow\infty}{\lim}\mu\left(  E_{n}\right)
=\mu\left( \underset{n=1}{\overset{\infty}{ {\displaystyle\bigcup}
}}E_{n}\right)  $ (continuity from below);

\item $\left\{  E_{n}\right\}  \subset\Sigma,E_{1}\supset E_{2}\supset
...,\mu\left(  E_{1}\right)  <\infty,$ imply
$\underset{n\rightarrow\infty }{\lim}\mu\left(  E_{n}\right)
=\mu\left(  \underset{n=1}{\overset{\infty}{ {\displaystyle\bigcap}
}}E_{n}\right)  $ (continuity from above).
\end{enumerate}

If $f$ is a non-negative real-valued function defined on $
X ,$ we denote the set
\[ \left\{  x\in X :f\left(
x\right)  \geq\alpha\right\} =\left\{ f\geq\alpha\right\}
\]
by $F_{\alpha}$ for $\alpha\geq0$. Note that if $\alpha\leq\beta$
then $F_{\beta}\subset F_{\alpha}.$

Let $\left(X, \Sigma, \mu\right)$ be a fuzzy measure space, we denote $M^{+}$ the set of all non-negative measurable functions with respect to $\Sigma$

\begin{dfn}\cite{Sugeno,Wang}
Let 
$A \in X$, 
$f\in M^{+}$
The fuzzy integral of $f$ on $A$ with respect to $\mu$ which is denoted by $\left(  s\right)
{\displaystyle\int}
fd\mu,$ is defined by
\[
\left(  s\right)
{\displaystyle\int}
fd\mu=\underset{\alpha\geq0}{%
{\displaystyle\bigvee}
}\left[  \alpha\wedge\mu\left(  A\cap\left\{  f\geq\alpha\right\}
\right) \right]  .
\] 
When $A=\Sigma,$ the fuzzy integral may also be denoted by $\left(  s\right)
{\displaystyle\int}
fd\mu.$\\
Where $\vee$ and $\wedge$ denote the operations inf and sup on
$[0,\infty)$, respectively.
\end{dfn}

The following properties of the Sugeno integral are well known and
can be found in.
\begin{prp}
Let $\left(X, \Sigma, \mu\right)$ be a fuzzy measure space, $A\in\Sigma$ and $f, g\in M^{+}$ 

\begin{enumerate}
\item $\left(  s\right)
{\displaystyle\int\limits_{A}}
fd\mu\leq\mu\left(  A\right)  $;

\item $\left(  s\right)
{\displaystyle\int\limits_{A}}
kd\mu=k\wedge\mu\left(  A\right),$ k non-negative constant;

\item If $f\leq g$ on $A$ then $\left(  s\right)
{\displaystyle\int\limits_{A}}
fd\mu\leq\left(  s\right)
{\displaystyle\int\limits_{A}}
gd\mu$;

\item $\mu\left(  A\cap\left\{  f\geq\alpha\right\}  \right)  \geq
\alpha\Rightarrow\left(  s\right)
{\displaystyle\int\limits_{A}}
fd\mu\geq\alpha$;

\item $\mu\left(  A\cap\left\{  f\geq\alpha\right\}  \right)  \leq
\alpha\Rightarrow\left(  s\right)
{\displaystyle\int\limits_{A}}
fd\mu\leq\alpha$;

\item $\left(  s\right)
{\displaystyle\int\limits_{A}}
fd\mu<\alpha\Leftrightarrow$ there exists $\gamma<\alpha$ such that
$\mu\left(  A\cap\left\{  f\geq\gamma\right\}  \right)  <\alpha$;

\item $\left(  s\right)
{\displaystyle\int\limits_{A}}
fd\mu>\alpha\Leftrightarrow$ there exists $\gamma>\alpha$ such that
$\mu\left(  A\cap\left\{  f\geq\gamma\right\}  \right)  >\alpha.$
\end{enumerate}
\end{prp}

\begin{rem}
Consider the distribution function F associated to $f$ on $A$, that
is, $F\left(  \alpha\right) =\mu\left( A\cap\left\{  f\geq
\alpha\right\}  \right)  .$ Then, due to (4) and (5) of Preposition
2.1, we have that
\[
F\left(  \alpha\right)  =\alpha\Rightarrow\left(
s\right)
{\displaystyle\int}
fd\mu=\alpha.
\]
Thus, from a numerical point of view, the fuzzy integral can be
calculated solving the equation $F\left(  \alpha\right)  =\alpha.$
\end{rem}

\cite{Caballero_4}, J.Caballero, K. Sadarangani proved with the help of certain
examples that the classical Hermite-Hadamard inequalities for fuzzy
integrals.

\begin{dfn}\cite{Dragomir_1} Let $I$ be an interval of real numbers.
A function $f:I\rightarrow(0,\infty)$ is said to be log-convex or
multiplicaticely convex if $\log\left(  f\right)  $ is convex, or,
equivalently, if for all $x,y\in I$ and $t\in\lbrack0,1]$ one has
the inequality:
\[
f\left(  tx+\left(  1-t\right)  y\right) \leq\left\vert f\left(
x\right)  \right\vert ^{t}\left\vert f\left( y\right) \right\vert
^{1-t}.
\]
\end{dfn}

We note that if $f$ and $g$ are convex functions and $g$ is
monotonic nondecrasing, then $g\circ f$ is convex. Moreover, since
$f=\exp\left( \log\left(  f\right)  \right)  ,$ it follows that a
log-convex function is convex, but the converse is not true.

\section{Hermite-Hadamard Type Inequality for Preinvex Functions via Sugeno Integrals}

The following inequality is well known in the literature as the
Hermite-Hadamard inequality
\[
f\left(  \frac{a+b}{2}\right)  \leq\frac{1}%
{b-a}\underset{a}{\overset{b}{\int}}f\left(  x\right)  dx\leq
\frac{f\left( a\right)  +f\left(  b\right)  }{2}
\]
where $f:I\rightarrow%
\mathbb{R}
$ is a convex function on the interval $I$ and $a,b\in I$ with
$a<b.$

In \cite{Dragomir_1}, S.S. Dragomir extended this classic result for
log-convex functions as follows:
\[
f\left(  \frac{a+b}{2}\right) \leq\frac{1}%
{b-a}\underset{a}{\overset{b}{\int}}f\left(  x\right)  dx\leq
L\left( f\left(  a\right)  ,f\left(  b\right)  \right),
\]
where $L\left(  p,q\right)  :=\dfrac{p-q}{\ln p-\ln q}\left(  p\neq
q\right) $ is the logaritmic mean of the positive real numbers $p,q$
(for $p=q,$ we put $L\left(  p,p\right)  =p$).

In this paper, we prove using Sugeno integral another refinement of
the Hermite-Hadamard type inequality for log-convex functions. Some
applications for special means are also given.\\

\noindent \textbf{Example 3.1.} Consider $X=[0,1]$ and let $\mu$ be
the Lebesgue measure on $X$. If we take the function $f\left(
x\right) =e^{x+1}$, then $f\left( x\right)$ is a log-convex
function. To calculate the Sugeno integral related to this function,
let's consider the distribution function $F$ associated to $f$ to
$[0,1]$, by Remark 2.1, this is
\begin{eqnarray*}
F\left(  \alpha\right)  &=&\mu\left( \lbrack0,1]\cap\left\{
f\geq\alpha\right\}  \right)  =\mu\left(  \lbrack0,1]\cap\left\{  e^{-x}%
\geq\alpha\right\}  \right) \\
&=&\mu\left( \lbrack0,1]\cap\left\{  x\leq -\ln\left(  \alpha\right)
\right\} \right)   =-\ln\left( \alpha\right).
\end{eqnarray*}
and we solve the equation
\[
-\ln\left(  \alpha\right)  =\alpha.\] It is easily proved that the
solutions of the last equation is 0.5672 with using bisection method
of numerical analysis, and, Remark 2.1, we get
\[\bigskip\qquad\left(  s\right)
\overset{1}{\underset{0}{\int}}fd\mu=\left( s\right)
\overset{1}{\underset{0}{\int}}e^{x+1}d\mu=0.5672 .\]
On the other
hand,
\[f\left(  \frac{0+1}{2}\right)  =f\left(  \frac{1}%
{2}\right)  =e^{-\frac{1}{2}}=0.6065 .\]

This proves that the left part of Hermite-Hadamard inequality is
not satisfied in the fuzzy context.\\

\noindent \textbf{Example 3.2.} Consider $X=[0,1]$ and let $\mu$ be
the Lebesque measure on $X$. Then for the log-convex function
$f\left( x\right) =e^{-\cos(x)-1}$ and using a similar argument that in
Example 3.1, we can get
\[
\left(  s\right)  \overset{1}{\underset{0}{\int}}fd\mu
=\left(  s\right) \overset{1}{\underset{0}{  \int}}\left(
e^{-\cos(x)-1}\right) d\mu=0.1852 \]
On the other hand,
\[L\left(  f\left(  0\right)  ,f\left(  1\right) \right)
=\dfrac{f\left(  0\right)  -f\left(  1\right)  }{\ln f\left(
0\right)  -\ln f\left(  1\right)  }=0.1718 \]
and this proves that
right-hand side of Hermite-Hadamard inequality is not satisfied for
fuzzy integrals.

The aim of the following theorem is to show a Hermite-Hadamard type
inequality for the Sugeno integral.

\begin{thm}
Let $g:[0,1]\rightarrow\lbrack0,\infty)$ be a
log-covex function such that $g\left(  0\right)  <g\left(  1\right)
$ and $\mu$ the
Lebesque measure on $%
\mathbb{R}
.$
Then
\[
\left(  s\right)  \overset{1}{\underset{0}{%
{\displaystyle\int}
}}gd\mu\leq\min\left\{  \alpha ,1\right\},
\]
where $\alpha=1-t$, $t$ satisfies the following equation
\[
\left[ g\left(0\right)\right]^{1-t}.
\left[ g\left(1\right)\right]^{t}+t-1=0
\]

\end{thm}

\noindent \textbf{Proof.} As a $g$ is a log-convex function, for
$x\in\lbrack0,1]$

\[g\left(  x\right)  =g\left(  \left(  1-x\right)
.0+x.1\right) \leq\left[  g\left(  0\right)  \right]  ^{1-x}.\left[
g\left(  1\right) \right]  ^{x}=h\left(  x\right)\] hence, by (3) of
Proposition 2.1, we have that
\begin{eqnarray*}
\left(  s\right)  \overset{1}{\underset{0}{%
{\displaystyle\int}
}}gd\mu&\leq&\ \left(  s\right) \overset{1}{\underset{0}{
{\displaystyle\int}
}}g\left(  \left(  1-x\right)  .0+x.1\right)  d\mu\\ \\
&\leq&\left(  s\right)  \overset{1}{\underset{0}{%
{\displaystyle\int}
}}\left[  g\left(  0\right)  \right]  ^{1-x}.\left[  g\left(
1\right)
\right]  ^{x}d\mu=\left(  s\right)  \overset{1}{\underset{0}{%
{\displaystyle\int}
}}h\left(  x\right)  d\mu.
\end{eqnarray*}
In order to calculate the integral in the right-hand part of the
last inequality, we consider the distribution function $F$ given by

\begin{eqnarray*}
F\left(  \alpha\right)  &=& \mu\left( \lbrack0,1]\cap\left\{
h\geq\alpha\right\}  \right)  =\mu\left( \lbrack0,1]\cap\left\{
\left[ g\left(  0\right)  \right] ^{1-x}.\left[  g\left(  1\right)
\right] ^{x}\geq\alpha\right\} \right)\\ \\
&=&\mu\left(
\lbrack0,1]\cap\left\{  x\geq\dfrac{\ln\left( \alpha\right)
-\ln\left( g\left(  0\right)  \right)  }{\ln\left( g\left(  1\right)
\right) -\ln\left(  g\left(  0\right)  \right)
}\right\}  \right)\\ \\
&=&1-\dfrac
{\ln\left(  \alpha\right)  -\ln\left(  g\left(  0\right)  \right)  }%
{\ln\left(  g\left(  1\right)  \right)  -\ln\left(  g\left( 0\right)
\right) },
\end{eqnarray*}
and the solution of the equation
\[1-\dfrac{\ln\left(  \alpha\right)  -\ln\left(
g\left( 0\right)  \right)  }{\ln\left(  g\left(  1\right)  \right)
-\ln\left( g\left(  0\right)  \right)  }=\alpha.\] Let $\alpha=1-t$, $t$ satisfies the following equation 
\[
\left[ g\left(0\right)\right]^{1-t}.
\left[ g\left(1\right)\right]^{t}+t-1=0.
\]
By (1) of
Proposition 2.1, we get that
\[\left(  s\right)  \overset{1}{\underset{0}{%
{\displaystyle\int}
}}h\left(  x\right)  d\mu\leq\mu\left(  \left[  0,1\right]  \right)
=1.\] By Remark 2.1, we have
\[\left(  s\right)  \overset{1}{\underset{0}{%
{\displaystyle\int}
}}gd\mu\leq\min\left\{  \alpha,1\right\}.\]
This completes is
proof.

\begin{rem}
In the case $g\left(  0\right)  =g\left( 1\right) $ in Theorem 3.1,
the function $h\left(  x\right)  $ is

\[h\left(  x\right)  =\left[  g\left(  0\right)
\right] ^{1-x}.\left[  g\left(  1\right)  \right]  ^{x}=g\left(
0\right) \]
and

\[\left(  s\right)  \overset{1}{\underset{0}{%
{\displaystyle\int}
}}gd\mu\leq\left(  s\right)  \overset{1}{\underset{0}{%
{\displaystyle\int}
}}h\left(  x\right)  d\mu=\left(  s\right)  \overset{1}{\underset{0}{%
{\displaystyle\int}
}}g\left(  0\right)  d\mu=g\left(  0\right)  \wedge1.\]
\end{rem}

\begin{thm}
Let $g:[0,1]\rightarrow\lbrack0,\infty)$ be a log-convex function
such that $g\left(  0\right)  >g\left(  1\right) $ and $\mu$ the
Lebesque measure on $%
\mathbb{R}
.$ Then

\[\left(  s\right)  \overset{1}{\underset{0}{%
{\displaystyle\int}
}}gd\mu\leq\min\left\{  \alpha,1\right\} \] where $\alpha$
is root of the equation
\[\dfrac{\ln\left(  \alpha\right)  -\ln\left( g\left(
0\right)  \right)  }{\ln\left(  g\left(  1\right)  \right)
-\ln\left(  g\left(  0\right)  \right)  }=\alpha.\]Let $\alpha=1-t$, $t$ satisfies the following equation 
\[
\left[ g\left(0\right)\right]^{t}.
\left[ g\left(1\right)\right]^{1-t}+t-1=0.
\] 
\end{thm}

\noindent \textbf{Proof.} Similarly, using the method in Theorem
3.1, we have

\begin{eqnarray*}
F\left(  \alpha\right)&=&\mu\left( \lbrack0,1]\cap\left\{
g\geq\alpha\right\}  \right)  \\
&=&\mu\left( \lbrack0,1]\cap\left\{ x \leq
\frac{\ln\left(  \alpha\right)  -\ln\left(  g\left(  0\right)  \right)  }%
{\ln\left(  g\left(  1\right)  \right)  -\ln\left(  g\left(
0\right) \right)  }\right\}  \right)\\
&=&\dfrac{\ln\left( \alpha\right)  -\ln\left(  g\left(  0\right)
\right)  }{\ln\left(  g\left( 1\right)  \right)  -\ln\left(  g\left(
0\right)  \right)  },
\end{eqnarray*}
and the solution of the equation
\[\dfrac{\ln\left(  \alpha\right)  -\ln\left(
g\left( 0\right)  \right)  }{\ln\left(  g\left(  1\right)  \right)
-\ln\left( g\left(  0\right)  \right)  }=\alpha,\] where $\alpha$ satisfies the following equation
\[
\left[ g\left(0\right)\right]^{1-\alpha}.
\left[ g\left(1\right)\right]^{\alpha}-\alpha=0.
\]
The proof of the rest part is similar, so we omit it. \\

\noindent \textbf{Example 3.3.} Consider $f\left(  x\right)
=e^{x^{2}-1}$ on $\left[ 0,1\right]  .$ Obviously, this function is
non-negative, non-decreasing and log-convex on the interval $\left[
0,1\right] $. Moreover, $f\left( 0\right)  =e^{-1}=\frac{1}{e}$ and $f\left(
1\right) =1>0.$ Calculating the fuzzy integral, we have
\[1-\dfrac{\ln\left( \alpha\right)  -\ln\left(
f\left(  0\right) \right)  }{\ln\left( f\left(  1\right)  \right)
-\ln\left(  f\left( 0\right)  \right) }=\alpha.\] Then, solving by
bisection method of numerical analysis, the approximately solution
$\alpha=0.5672.$ By Theorem 3.1, we have
\[\left(  s\right)  \overset{1}{\underset{0}{%
{\displaystyle\int}
}}fd\mu\leq\min\left\{  \alpha,1\right\}  =0.5672.\]
Also $t$ is the root of the $\alpha=1-t$ equation, satisfies the following equation
\[
e^{t-1}+t-1=0.
\] 
\noindent \textbf{Example 3.4.} Consider the log-convex function
$f\left( x\right) =e^{-\sin(x)}$, for $x\in\left[  0,1\right]. $ Then
$f\left( 0\right)  =1$ and $f\left( 1\right)  =0.4311,$ and we have
\[\dfrac{\ln\left(  \alpha\right)
-\ln\left( f\left( 0\right) \right)  }{\ln\left(  f\left(  1\right)
\right) -\ln\left( f\left( 0\right)  \right)  }=\alpha\] which gives
by solving by bisection method of numerical analysis, the
approximately solution $\alpha=0.6024,$ satisfies under the equation 
\[
\ln(\alpha)+\sin(1)*\alpha=0.
\]
By Theorem 3.2, we
have estimate:
\[\left(  s\right)  \overset{1}{\underset{0}{%
{\displaystyle\int}
}}e^{-\sin(x)}d\mu\leq\min\left\{  \alpha,1\right\}  =\alpha.\]

\begin{thm}Let $g:\left[  a,b\right]
\rightarrow\lbrack0,\infty)$ be
a log-convex function and $\mu$ the Lebesque measure on $%
\mathbb{R}
.$ Then\\

\indent\text{(i) }If $g\left(  a\right)  <g\left(  b\right) $, then

\[\left(  s\right)  \overset{b}{\underset{a}{%
{\displaystyle\int}
}}gd\mu\leq\min\left\{  \alpha_{1},b-a\right\}\]
where $\alpha_{1}$
is root of the equation
\[b-\dfrac{\left( b-a\right)  .\ln\left(
\alpha\right) -b.\ln\left(  g\left( a\right)  \right)  +a.\ln\left(
g\left( b\right)  \right) }{\ln\left(  g\left(  b\right)  \right)
-\ln\left( g\left( a\right)  \right)  }=\alpha.\]

\indent\text{(ii)} If $g\left(  a\right)  =g\left(  b\right)  $,
then
\[\left(  s\right)  \overset{b}{\underset{a}{%
{\displaystyle\int}
}}gd\mu\leq\min\left\{  g\left(  a\right)  ,b-a\right\}.  \]

\indent\text{(iii) }If $g\left(  a\right)  >g\left(  b\right) $,
then
\[\left(  s\right)  \overset{b}{\underset{a}{%
{\displaystyle\int}
}}gd\mu\leq\min\left\{  \alpha_{2},b-a\right\}, \]
where
$\alpha_{2}$ is root of the equation
\[\dfrac{\left(  b-a\right)
.\ln\left( \alpha\right) -b.\ln\left(  g\left(  a\right)  \right)
+a.\ln\left( g\left( b\right) \right)  }{\ln\left(  g\left( b\right)
\right) -\ln\left( g\left( a\right)  \right) }-a=\alpha. \]
\end{thm}

\noindent \textbf{Proof.} We will prove (i) and other two cases are
similar. Note that as $g$ is a log-convex function then for
$x\in\left[ 0,1\right]  $ we have

\[\bigskip\qquad g\left(  x\right)  =g\left(  \left(
1-\frac{x-a}{b-a}\right) .a+\frac{x-a}{b-a}.b\right)  \leq\left(
g\left(  a\right)  \right)
^{\frac{b-x}{b-a}}.\left(  g\left(  b\right)  \right)  ^{\frac{x-a}{b-a}%
}=h\left(  x\right).\] By (3) of Proposition 2.1,
\[\qquad\left(  s\right)  \overset{b}{\underset{a}{%
{\displaystyle\int}
}}gd\mu\leq\left(  s\right)  \overset{b}{\underset{a}{%
{\displaystyle\int}
}}\left(  g\left(  a\right)  \right)  ^{\frac{b-x}{b-a}}.\left(
g\left( b\right)  \right)  ^{\frac{x-a}{b-a}}d\mu=\left(  s\right)
\overset
{b}{\underset{a}{%
{\displaystyle\int}
}}h\left(  x\right)  d\mu.\]
Now, we consider the distribution
function $F$ given by
\begin{eqnarray*}
F\left(  \alpha\right)  &=&\mu\left(  \left[  a,b\right] \cap\left\{
h\geq\alpha\right\}  \right)  =\mu\left(  \left[ a,b\right]
\cap\left\{ \left(  g\left(  a\right)  \right)
^{\frac{b-x}{b-a}}.\left(  g\left( b\right)  \right)
^{\frac{x-a}{b-a}}\geq\alpha\right\}  \right)\\ \\
&=&\mu\left(
\left[ a,b\right] \cap\left\{  x\geq\dfrac{\left( b-a\right)
.\ln\left( \alpha\right) -b.\ln\left(  g\left(  a\right) \right)
+a.\ln\left( g\left(  b\right) \right)  }{\ln\left( g\left( b\right)
\right) -\ln\left(  g\left( a\right)  \right)  }\right\}
\right)\\ \\
&=&b-\dfrac{\left( b-a\right) .\ln\left(  \alpha\right)
-b.\ln\left( g\left( a\right) \right) +a.\ln\left(  g\left( b\right)
\right) }{\ln\left( g\left( b\right) \right)  -\ln\left( g\left(
a\right) \right)  },\\
\end{eqnarray*}
and the root is $\alpha_{1}$ which is the solution of the equation

\[b-\dfrac{\left(  b-a\right)  .\ln\left(  \alpha\right)
-b.\ln\left(  g\left(  a\right)  \right)  +a.\ln\left(  g\left(
b\right) \right)  }{\ln\left(  g\left(  b\right)  \right) -\ln\left(
g\left( a\right)  \right) }=\alpha.\]

Then by (1) of Proposition 2.1 and Remark 2.1, we have

\[\left(  s\right)  \overset{b}{\underset{a}{%
{\displaystyle\int}
}}gd\mu\leq\left(  s\right)  \overset{b}{\underset{a}{%
{\displaystyle\int}
}}h\left(  x\right)  d\mu=\min\left\{  \alpha_{1}, b-a\right\}  .\]\\

\noindent \textbf{Example 3.5.} Consider $f\left(  x\right)
=e^{-\sin\left(2x\right) }$ be a function defined on $\left[
\frac{\pi}{4},\frac{\pi}{2} \right] $. This function is non-decreasing and
log-convex because $\log(\exp(-sin(2x)))=-\sin(2x) $ function is convex and $f\left(  x\right)
=e^{-\sin\left( 2x\right) }$ is non-negative. As $f\left( \frac{\pi}{4}\right)  =0.3679$ and $f\left(
\frac{\pi}{2}\right) =1$ and $f\left( \frac{\pi}{4}\right) <$ $f\left(
\frac{\pi}{2}\right)  $, by (a) of Theorem 3.3 we can get the
following estimate:

\[\left(  s\right)  \overset{\pi/2}{\underset{\pi/4}{%
{\displaystyle\int}
}}\left(  e^{-\sin\left( 2x\right)  }\right)  d\mu\leq\min\left\{
\alpha_{1},\frac{\pi}{4}\right\}  \] where $\alpha_{1}$ is root
which is the equation

\[\frac{\pi}{4}-\dfrac{\left(  \frac{\pi}%
{4}\right)  .\ln\left(  \alpha\right)  -\frac{\pi}{2}.\ln\left(
g\left(
\frac{\pi}{4}\right)  \right)  +\frac{\pi}{4}.\ln\left(  g\left(  \frac{\pi}{2}\right)  \right)  }%
{\ln\left(  g\left(  \frac{\pi}{2}\right)  \right)  -\ln\left(
g\left(\frac{\pi}{4}\right)  \right)  }=\alpha.\]

This equation have been solved by matlab program and the root is
$\alpha_{1}=0.5175.$ Definitively Sugeno integral:

\[\left(  s\right)  \overset{\pi/4}{\underset{\pi/4}{%
{\displaystyle\int}
}}\left(  e^{-\sin\left(  2x\right)  }\right)  d\mu\leq\min\left\{
\alpha_{1},\frac{\pi}{4}\right\}  =\alpha_{1}=0.5175.\]

\section{Conclusion}

In this paper, we have researched the classical Hermite-Hadamard
inequality for Sugeno integral based on log-convex function. For
further investigations we will continue to study Hermite-Hadamard
and other integral inequalities for several fuzzy integrals
based on log-convex function.

\section{Acknowledgments}

The authors would like to thank to Assist. Prof. Imdat ISCAN, Department of Mathematics Giresun University, owing to his ideas and would like to thank to the referees for reading this work carefully, providing valuable suggestions and comments.

\end{document}